\documentclass[12pt,reqno]{amsart}

\usepackage{amssymb, amsmath, amsfonts, latexsym}
\usepackage{enumerate}
\usepackage{color}

\newtheorem{thm}{Theorem}[section]
\newtheorem{cor}[thm]{Corollary}
\newtheorem{lem}[thm]{Lemma}
\newtheorem{prop}[thm]{Proposition}

\newtheorem{remarks}[thm]{Remark}

\theoremstyle{definition}
\newtheorem{defn}{Definition}[section]

\numberwithin{equation}{section} \theoremstyle{remark}

\newcommand{\rr}{\mathbb{R}}
\newcommand{\R}{\mathbb{R}}

\newcommand{\s}{\mathbb{S}}

\def\inte{{\mathrm{{\rm Int}}}}
\def\AA{\mathcal A}

\def\CC{\mathcal C}
\def\DD{\mathcal D}

\def\EE{\mathcal E}

\def\LL{\mathcal L}
\def\MM{\mathcal M}

\def\VV{\mathcal V}

\def\R{\mathbb R}

\def\vep{\varepsilon}
\def\<{\langle}
\def\>{\rangle}

\def\beq{\begin{equation}}
\def\nneq{\end{equation}}

\def\bdef{\begin{defn}}
\def\ndef{\end{defn}}

\def\bthm{\begin{thm}}
\def\nthm{\end{thm}}

\def\bprop{\begin{prop}}
\def\nprop{\end{prop}}

\def\brmk{\begin{remarks}}
\def\nrmk{\end{remarks}}

\def\bexa{\begin{exa}}
\def\nexa{\end{exa}}

\def\blem{\begin{lem}}
\def\nlem{\end{lem}}

\def\bcor{\begin{cor}}
\def\ncor{\end{cor}}

\def\bexe{\begin{exe}}
\def\nexe{\end{exe}}

\def\bprf{\begin{proof}}
\def\nprf{\end{proof}}

\def\bdes{\begin{description}}
\def\ndes{\end{description}}

\def\Var{{\rm Var}}

\begin{document}

\title[Spectral gap]{\bf Spectral gap for spherically symmetric log-concave probability measures, and beyond}

\author{Michel Bonnefont}
\address{Michel BONNEFONT \\ Institut de Math\'ematiques de Bordeaux, Universit\'e  de Bordeaux, 351 cours de la Lib\'eration,
33405 Talence, France.
}
\thanks{MB is partially supported by the French ANR-12-BS01-0013-02 HAB project}
\email{michel.bonnefont@math.u-bordeaux.fr}
\urladdr{http://www.math.u-bordeaux.fr/~mibonnef/}

\author{Ald\'eric Joulin}
\address{Ald\'eric JOULIN (corresponding author) \\ Universit\'e de Toulouse,
Institut National des Sciences Appliqu\'ees,
Institut de Math\'ematiques de Toulouse,
F-31077 Toulouse,
France.
}
\thanks{AJ is partially supported by the French ANR-2011-BS01-007-01 GEMECOD and ANR-12-BS01-0019 STAB projects}
\email{ajoulin@insa-toulouse.fr}
\urladdr{http://www.math.univ-toulouse.fr/~joulin/}

\author{Yutao Ma}
\address{Yutao MA\\ School of Mathematical Sciences $\&$ Lab. Math. Com. Sys., Beijing Normal University, 100875 Beijing, China.}
\thanks{YM is partially supported by NSFC 11101040, 11371283, 985 Projects and the Fundamental Research Funds for the Central Universities.}
\email{mayt@bnu.edu.cn}
\urladdr{http://math.bnu.edu.cn/~mayt/}

\keywords{Spectral gap; Diffusion operator; Poincar\'e-type inequalities; Log-concave probability measure.}

\subjclass[2010]{60J60, 39B62, 37A30, 58J50.}

\maketitle

\begin{abstract}
Let $\mu$ be a probability measure on $\rr^n$ ($n \geq 2$) with Lebesgue density proportional to $e^{-V (\Vert x\Vert )}$, where $V : \rr_+ \to \rr$ is a smooth convex potential. We show that the associated spectral gap in $L^2 (\mu)$ lies between $(n-1) / \int_{\rr^n} \Vert x\Vert ^2 \mu(dx)$ and $n / \int_{\rr^n} \Vert x\Vert ^2 \mu(dx)$, improving a well-known two-sided estimate due to Bobkov. Our Markovian approach is remarkably simple and is sufficiently robust to be extended beyond the log-concave case, at the price of potentially modifying the underlying dynamics in the energy, leading to weighted Poincar\'e inequalities. All our results are illustrated by some classical and less classical examples.
\end{abstract}

\section{Introduction}

Given a complete Riemannian manifold endowed with a smooth probability measure $\mu$, we associate a canonical non-negative diffusion operator $- \LL_\mu$ that is self-adjoint in $L^2 (\mu)$. Then a natural question arises: does the Markovian operator $-\LL _\mu$ have a gap in the bottom of its spectrum ? If yes, can we estimate the size of this gap, in particular with respect to the dimension of the manifold ? Since the pioneer work of Lichnerowicz at the end of the fifties, these questions have attracted the attention of many mathematicians in the last few decades. If $\lambda _1 (- \LL _\mu) $ denotes the size of this gap, then it can be alternatively formulated as the optimal constant in the Poincar\'e inequality. Such an object, which is a key ingredient in the theory of functional inequalities, controls the variance of functions in terms of the $L^2$-norm of their gradient. From analytical and probabilistic perspectives, the existence of a spectral gap is of crucial importance since on the one hand it governs the convergence to equilibrium of the underlying Markov process and on the other hand it has fundamental consequences in the so-called concentration of measure phenomenon but also in isoperimetric problems and their deep connections to Ricci curvature in geometry. The interested reader is referred to the set of notes of Ledoux \cite{ledoux_berlin,Ledoux} and also to the recent monograph \cite{BGL} for a nice introduction to the topic, with historical references and precise credit for this large body of work.

In the case of probability measures $\mu$ on the Euclidean space $(\rr^n , \Vert \cdot \Vert )$ with Lebesgue density proportional to $e^{-V}$, where $V$ is some smooth potential, a natural Markovian operator is given for all smooth enough functions $f$ by
$$
\LL _\mu f = \Delta f - \nabla V \cdot \nabla f .
$$
Using for instance the famous Bakry-\'Emery criterion \cite{bakry_emery} on the basis of the so-called $\Gamma_2$-calculus, a positive lower bound $K$ on the Hessian of $V$ ensures the existence of the spectral gap and moreover it is at least $K$. However, if $V$ is assumed to be only convex or if it exhibits some convexity only at infinity (allowing concavity in a localized region of the space), then there might be a spectral gap whereas the Bakry-\'Emery criterion fails. In the last thirty years, there has been many rooms of improvements of such a result and in many different directions. Recently, Bobkov \cite{bobkov_spherical} obtained a nice estimate on the spectral gap for a class of log-concave probability measures. More precisely, such measures $\mu$ are spherically symmetric and log-concave, meaning that the potential $V$ is radial and convex: it can be written as $V(x) = V(\Vert x \Vert)$ with $V : (0,+\infty) \to (0,+\infty)$ a convex function (by abuse of notation, we shall use all along the paper the same notation for the functions defined on $\rr^n$ and their radial version defined on $(0,+\infty)$). Then the operator above rewrites as
\[
\LL _\mu f( x)  = \Delta f (x) - \frac{V' (\Vert x\Vert )}{\Vert x\Vert } \, x \cdot \nabla f( x) , \quad x\in \rr^n ,
\]
and then his famous theorem can be stated as follows.
\bthm[Bobkov's theorem \cite{bobkov_spherical}]
\label{theo:bobkov}
Let $\mu$ be a spherically symmetric log-concave probability measure on $\rr^n$ ($n \geq 2$). Then the spectral gap $\lambda_1 (-\LL_\mu)$ satisfies
$$
\frac{n}{13 \int_{\rr ^n} \Vert x\Vert ^2 \mu (dx)} \leq \lambda_1 (-\LL_\mu) \leq \frac{n}{\int_{\rr ^n} \Vert x\Vert ^2 \mu (dx)}  .
$$
\nthm

A consequence of Bobkov's theorem is that the class of spherically symmetric log-concave probability measures satisfies the KLS conjecture: if the covariance of $\mu$ is the identity, meaning that for all $\theta \in \rr ^n$,
$$
\int _{\rr ^n} \vert x \cdot \theta \vert ^2 \mu (dx) = \Vert \theta \Vert ^2,
$$
then the spectral gap $\lambda_1 (-\LL_\mu)$ is bounded from below by an universal constant, independent from the dimension $n$. Recall that the KLS conjecture was originally introduced in \cite{KLS} for the Cheeger inequality but it is actually equivalent to state it for the Poincar\'e inequality according to a result of Ledoux \cite{Ledoux}, generalized recently by Milman \cite{milman}. \vspace{0.1cm}

The next theorem is the main result of the paper and is an improvement of Bobkov's theorem. In particular it gives the exact asymptotics of the spectral gap as the dimension of the space goes to infinity.
\bthm[Bobkov's theorem revisited]
\label{theo:main}
Under the assumptions and notation of Bobkov's theorem, the spectral gap $\lambda_1 (-\LL_\mu)$ satisfies
$$
\frac{n-1}{\int_{\rr ^n} \Vert x\Vert ^2 \mu (dx)} \leq \lambda_1 (-\LL_\mu) \leq \frac{n}{\int_{\rr ^n} \Vert x\Vert ^2 \mu (dx)}  .
$$
\nthm
Certainly, such a result can be adapted without additional effort to the case of a measure $\mu$ whose support is a centered Euclidean ball provided Neumann's boundary conditions are assumed for the underlying dynamics (we will investigate such an example later). \vspace{0.1cm}

As we will see in the sequel, there are two main ingredients in the proof of Theorem \ref{theo:main} which reveals to be remarkably simple. The first one is a comparison between the spectral gap of the full dynamics $-\LL_\mu$ and that of its one-dimensional radial part. Such a comparison is conceivable since the probability measure $\mu$ is close to the product measure whose marginals are its spherical and radial parts, and because of the well-known tensorization property of the Poincar\'e inequality. The second main point of the proof of Theorem \ref{theo:main} focuses on a careful estimation of the spectral gap of the one-dimensional radial part, based on a simple, but somewhat very useful, result given recently by two of the three authors \cite{bj}. \vspace{0.1cm}

Before entering into more details, let us provide the organization of the paper. In Section \ref{sect:prelim}, we briefly introduce some basic material on Markovian diffusion operators, their spectral gap and the links with the Poincar\'e inequality. Section \ref{sect:proof} is devoted to the proof of our main result, Theorem \ref{theo:main}, which is illustrated in Section \ref{sect:ex} by some examples such as probability measures $\mu$ with Lebesgue densities proportional to $e^{- \Vert x\Vert ^\alpha /\alpha}$, where $\alpha$ is a parameter greater than or equal to 1 to ensure the log-concavity of the measure $\mu$. For $\alpha = 1$ it corresponds to the Euclidean version of the exponential distribution and for $\alpha = 2$ to the classical Gaussian setting, whereas we obtain the uniform measure on the Euclidean unit ball as $\alpha \to +\infty$. Finally, we extend in Section \ref{sect:beyond} this technology to non-necessarily log-concave distributions by modifying in a convenient way the original dynamics, leading to weighted Poincar\'e inequalities. To observe the relevance of our approach beyond the log-concave case, we investigate in detail the example of heavy-tailed probability measures such as the generalized Cauchy distributions and we finish this work by revisiting the Gaussian case.

\section{Preliminaries}
\label{sect:prelim}

\subsection{Generalities on Markovian diffusion operators}
Let $(\MM,g)$ be a smooth complete connected $n$-dimensional Riemannian manifold with metric $g$. Denote $\mu$ the probability measure with density with respect to the volume measure of $(\MM,g)$ proportional to $e^{-V}$, where $V$ a smooth function on $(\MM,g)$. Let $\Delta_\MM$ be the Laplace-Beltrami operator on $(\MM,g)$ and denote $\LL _\mu$ the Markovian diffusion operator initially defined on the space $\CC ^\infty(\MM)$ of infinitely differentiable real-valued functions on $(\MM,g)$ by
$$
\LL_\mu f = \Delta _\MM f - \nabla _\MM V \cdot \nabla _\MM f ,
$$
where $\nabla_\MM $ stands for the Riemannian gradient and the underlying scalar product $\cdot$ is given with respect to the metric $g$ (denote $\vert \cdot \vert $ the associated norm). Let $\CC_0 ^\infty(\MM)$ be the subspace of $\CC ^\infty(\MM)$ consisting of compactly supported functions. The operator $\LL_\mu$ is symmetric on $\CC_0 ^\infty(\MM)$ with respect to the measure $\mu$, i.e. for every $f_1,f_2 \in \CC _0 ^\infty(\MM)$,
\[
\EE _\mu (f_1,f_2) := \int_\MM f_1 (-\LL_\mu f_2) d\mu = \int_\MM (-\LL _\mu f_1) f_2 d\mu  = \int_\MM \nabla _\MM f_1 \cdot \nabla _\MM f_2 d\mu,
\]
hence $\LL_\mu$ is non-positive on $\CC_0 ^\infty(\MM)$. Since the manifold is complete, a straightforward adaptation of the argument in \cite{strichartz} (see for example \cite{grigoryan} and the references therein) shows that $(\LL_\mu , \CC_0 ^\infty(\MM))$ is essentially self-adjoint in $L^2 (\mu)$, that is, it admits a unique self-adjoint extension (still denoted $\LL_\mu$) with domain $\DD (\LL_\mu) \subset L^2(\mu)$ in which the space $\CC_0 ^\infty(\MM)$ is dense for the norm
$$
\Vert f\Vert _{\LL_\mu} := \left( \Vert f \Vert ^2 _{L^2 (\mu)} + \Vert \LL_\mu f \Vert ^2 _{L^2 (\mu)} \right) ^{1/2} .
$$
In other words the space $\CC_0 ^\infty(\MM)$ is a core of the domain $\DD (\LL_\mu)$. In general, it is not easy to describe its domain, however it is easier to describe that of the associated Dirichlet form. The closure $(\EE _\mu , \DD (\EE _\mu))$ of the bilinear form $(\EE _\mu , \CC_0 ^\infty(\MM))$ is a Dirichlet form on $L^2(\mu)$ and we have the dense inclusion $\DD (\LL_\mu) \subset \DD (\EE _\mu)$ for the norm
$$
\Vert f\Vert _{\EE _\mu} := \left( \Vert f \Vert ^2 _{L^2 (\mu)} + \EE _\mu (f,f) \right) ^{1/2} ,
$$
where the domain $\DD (\EE _\mu)$ is the Sobolev space
\[
W^{1,2}(\mu) = \left\{ f\in L^2(\mu) : \vert \nabla _\MM f \vert \in  L^2(\mu) \right\} .
\]
If the manifold $(\MM,g)$ admits a smooth boundary, one can still consider the operator $\LL_\mu$ provided Neumann boundary conditions are assumed (in this case, the operator has to be defined initially on the space of smooth functions with vanishing normal derivative at the boundary). We refer to \cite{davies} for more details.

\subsection{Spectral gap} Let us consider the spectrum in $L^2(\mu)$ of the operator $-\LL_\mu$. First, it is a non-negative self-adjoint operator and thus its spectrum is included in $[0,+\infty)$. Since $\mu$ is assumed to be a probability measure, the constant functions belong to the domain $\DD (\LL_\mu)$ and the identity $\LL_\mu 1 = 0$ means that 0 is an eigenvalue, the constant functions being the associated eigenfunctions (by eigenvalue and eigenfunction, we mean the numbers $\lambda$ and the non-null functions $f\in \DD (\LL_\mu)$ satisfying $-\LL_\mu f = \lambda f$). In particular, the ergodicity property holds, i.e., the only functions $f$ in the domain $\DD (\LL_\mu)$ satisfying $\LL_\mu f = 0$ are the constant functions. Sometimes, it may happen that the operator $- \LL_\mu$ exhibits a gap in the bottom of its spectrum. More precisely, the spectral gap of the operator $- \LL_\mu$, usually denoted $\lambda_1 (- \LL _\mu)$, is the largest real $\lambda$ such that its spectrum lies in $\{0\} \cup [\lambda,+\infty)$. In other words, it is characterized by the variational identity
\begin{equation}
\label{eq:variational}
\lambda_1 (- \LL _\mu) = \inf \left \{ \int_\MM f(-\LL_\mu f) d\mu : f \in \DD (\EE_\mu) , \, f \perp  \mathbf 1, \, \int_\MM f^2 d\mu = 1 \right \} ,
\end{equation}
where $f \perp \mathbf 1$ means that $f$ is orthogonal to the constants, that is, $\int_\MM f d\mu =0$. Among the potential cases of interest we have in mind, one of the most famous example is the Laplace-Beltrami operator on the unit sphere $\s_{n-1} \subset \rr ^n$ ($n \geq 2$), endowed with the normalized Hausdorff measure $\sigma_{n-1}$ ($V\equiv 0$ in this case) and for which we have
$$
\lambda_1 (-\Delta _{\s ^{n-1}}) = n-1 .
$$
Writing \eqref{eq:variational} slightly differently, the spectral gap is characterized by the best (i.e., largest) constant $\lambda$ in the following Poincar\'e inequality,
\begin{equation}
\label{eq:poincare}
\lambda \Var _\mu (f) \leq \int _\MM  f (-\LL _\mu f) d\mu = \int_\MM \vert \nabla _\MM f \vert ^2 d\mu ,
\end{equation}
for all $f\in \DD(\EE_\mu)$ or, equivalently, for all $f\in \CC_0^\infty(\MM)$ by self-adjointness. Above $\Var _\mu (f)$ stands for the variance of $f$ under $\mu$, that is,
$$
\Var _\mu (f) := \int_\MM f^2 d\mu - \left( \int_\MM f d\mu \right) ^2 .
$$
As mentioned in the Introduction, such an inequality reveals to be very useful for probabilistic issues such as the $L^2$-convergence to equilibrium for the underlying Markov process, but also in analysis through measure concentration, isoperimetry and their links to geometry via the key role of the Ricci curvature. 
\section{Preparation and proof of Theorem \ref{theo:main}}
\label{sect:proof} This part is devoted to the preparation and the proof of our main contribution Theorem \ref{theo:main} stated in the Introduction. Recall that the underlying Markovian operator of interest acts on functions $f \in \CC_0^\infty(\R^n)$ as
\begin{equation}
\label{eq:operator_lmu}
\LL _\mu f( x): = \Delta f (x) - \frac{V' (\Vert x\Vert )}{\Vert x\Vert } \, x \cdot \nabla f( x) , \quad x\in \rr^n ,
\end{equation}
where $\Delta $ and $\nabla$ denote the Euclidean Laplacian and gradient, respectively. The potential $V$ is assumed to be convex and the symmetric and non-positive operator $(\LL_\mu , \CC_0^\infty(\R^n))$ with respect to the measure $\mu$ is essentially self-adjoint in $L^2(\mu)$. \vspace{0.1cm}

The proof of Theorem \ref{theo:main} is based on two complementary lemmas, which follow a specific strategy:
\begin{itemize}
\item[$\circ$] to understand the competition between the spectral gaps of the full dynamics and its one-dimensional radial part;
\item[$\circ$] to analyse and estimate the spectral gap of the radial part.
\end{itemize}

\subsection{Spectral gap estimates for Sturm-Liouville dynamics}
We begin by the second point and consider the one-dimensional case. Let $\nu$ be a one-dimensional probability measure whose support is some interval $I$, possibly finite. We assume that its Lebesgue density on $I$ is proportional to $e^{-U}$, where $U$ is a smooth convex potential on $I$. Then the log-concave measure $\nu$ is invariant for the following Sturm-Liouville dynamics
\begin{equation}
\label{eq:operator_lnu}
\LL _\nu f(r) : = f''(r) - U'(r) f'(r) , \quad r \in \inte (I) ,
\end{equation}
where $\inte (I)$ denotes the interior of $I$. Obviously, if $I$ admits a boundary then we assume Neumann's conditions. Denote $\lambda_1 (-\LL_\nu)$ the associated (Neumann) spectral gap corresponding to the optimal constant in the following Poincar\'e inequality: for any $f \in W^{1,2}(\nu)$,
\begin{equation}
\label{eq:poincare1d}
\lambda \Var _\nu (f) \leq \int _I  f (-\LL _\nu f) d\nu = \int_I ( f' )^2 d\nu .
\end{equation}
Among the various criteria available to estimate the spectral gap in dimension 1 (not only in the log-concave case), one has recently been put forward in \cite{bj} by two of the three authors and can be reformulated as follows: if $\CC^\infty (I)$ stands for the space of real-valued smooth functions on $I$, then we have
\begin{equation}
\label{e:chen}
\lambda_1 (-\LL_\nu) \geq \sup_{f\in \mathcal F} \inf_{x\in I}  V_f (x) ,
\end{equation}
where
\[
 V_f := - \frac{(\LL _\nu f)'} {f'} \quad \textrm{ with } \quad \mathcal F := \left\{ f\in \CC^\infty (I) : f' > 0 \mbox{ on } \inte (I) \right\} .
\]
Moreover the equality holds in \eqref{e:chen} if $\lambda_1 (-\LL_\nu)$ is an eigenvalue. Such a result already appeared via another approach in the work of Chen and Wang \cite{chen_wang} (see also \cite{chen} for an alternative form of this variational formula and \cite{djellout} which turns the form into the present one) and has been rediscovered in \cite{bj} through intertwining relations between gradients and semigroups. In general, it provides ``easy-to-verify" conditions ensuring the existence of a spectral gap for the dynamics. However, this formula can be difficult to use in order to obtain a good estimate since one has to guess what the best function $f$ is, and this is nothing but the eigenfunction (when it exists) associated to the spectral gap. \vspace{0.1cm}

To overcome this difficulty, at least in the log-concave case, the following result derived from \cite{bj} is well-adapted.
\blem[Bonnefont-Joulin]
\label{lemme:bj}
Assume that the second derivative of $U$ is positive on $\inte (I)$. Then the spectral gap $\lambda_1 (-\LL_\nu)$ of the operator $- \LL_\nu$ defined in \eqref{eq:operator_lnu} satisfies
\begin{equation}
\label{eq:bj_gen}
\lambda_1 (-\LL_\nu) \geq \frac{1}{\int_I \frac{1}{U''} d\nu } .
\end{equation}
In particular if $\nu$ has Lebesgue density on $\rr_+$ proportional to $r^{n-1} e^{- V(r)}$ ($n\geq 2$) and $V$ is a smooth convex function, then we have
\begin{equation}
\label{eq:bj_radial}
\lambda_1 (-\LL_\nu) \geq \frac{n-1}{\int_0 ^{+\infty} r^2 \nu (dr) } .
\end{equation}
\nlem
The first estimate above is an integrated version of the one-dimensional Bakry-\'Emery criterion. Moreover, it can be seen as a one-dimensional version of a nice result of Veysseire \cite{veysseire} given in the context of compact Riemannian manifolds, our second derivative $U''$ being replaced by the Ricci curvature lower bound. Although the inequality \eqref{eq:bj_gen} is sharp in some cases like for the standard Gaussian distribution associated to the Ornstein-Uhlenbeck process on $\rr$, it can sometimes be non-optimal. For instance if $1/U''$ is not $\nu$-integrable, then the latter lower bound does not provide any information on the spectral gap. However, it reveals to be convenient at least when $U''$ tends to 0 sufficiently slowly at the boundary of $I$ (at infinity if $I$ is unbounded), since the Bakry-\'Emery criterion is not available in this context (we would have $\lambda_1 (-\LL_\nu) \geq \inf U'' = 0$). \vspace{0.1cm}

The second inequality \eqref{eq:bj_radial} is the one we are interested in. The measure $\nu$, which belongs to the class of the so-called log-concave distributions of order $n$ according to Bobkov's definition \cite{bobkov_gaussian}, is the radial part of the spherically symmetric measure $\mu$ introduced in Theorem \ref{theo:main}. In particular, due to the presence of the important prefactor $r^{n-1}$, a log-concave distribution of order $n$ is ``more" log-concave than a classical log-concave measure. Hence we expect some interesting consequences of the general inequality \eqref{eq:bj_gen} applied to this class of log-concave measures, and this is indeed the case: setting
$$
U(r) := V(r) - (n-1) \log (r), \quad r >0,
$$
one has
$$
U''(r) = V''(r) + \frac{ n-1}{r^2} \geq \frac{ n-1}{r^2}  , \quad r >0,
$$
so that the inequality \eqref{eq:bj_radial} is a straightforward consequence of \eqref{eq:bj_gen}.

\subsection{Spectral comparison with the radial part}
As announced, we turn now to the comparison between the spectral gaps of the full dynamics and its one-dimensional radial part. We mention that the next result holds for general spherically symmetric probability measure (not necessarily log-concave) admitting a finite second moment.
\blem[Spectral comparison]
\label{lemme:gen-spherically}
Let $\mu$ be a spherically symmetric probability measure on $\R^n$ ($n \geq 2$) with Lebesgue density proportional to $e^{-V(\Vert x\Vert )}$ and let $ \nu$ be its radial part, that is, the one-dimensional probability measure on $(0,+\infty)$ whose Lebesgue density is  proportional to $r^{n-1} e^{- V(r)}$. Then the spectral gaps $\lambda_1 (-\LL_\mu)$ and $\lambda_1 (-\LL_\nu)$ of the operators $-\LL_\mu$ and $-\LL_\nu$ defined in \eqref{eq:operator_lmu} and \eqref{eq:operator_lnu} respectively, satisfy
\[
\min \left \{ \lambda_1 (-\LL_\nu) , \frac{n-1}{\int_{\rr ^n} \Vert x\Vert ^2 \mu (dx) }\right \}  \leq \lambda_1 (-\LL_\mu) \leq \min \left \{ \lambda_1 (-\LL_\nu) , \frac{n}{\int_{\rr ^n} \Vert x\Vert ^2 \mu (dx) }\right \}.
\]
\nlem
\bprf
The upper bound in the right-hand-side of the desired estimates is straightforward. Indeed, considering only radial functions, the variational characterization \eqref{eq:variational} of the spectral gap gives
\[
\lambda_1 (-\LL_\mu) \leq  \lambda_1 (-\LL_\nu),
\]
whereas the second part is obtained by choosing in the Poincar\'e inequality \eqref{eq:poincare} the linear function $f(x) = \sum_{k=1} ^n x_k$. \vspace{0.1cm}

Let us concentrate our attention on the challenging lower bound. The proof given below is a refinement of Bobkov's argument used to derive his famous theorem.
The principle is to benefit from the spherically invariance of $\mu$ and then to use the tensorization property of the Poincar\'e inequality
with respect to the radial and spherical measures. More precisely, the measure $\mu$ is the image measure
of the product measure $\nu \otimes \sigma_{n-1}$ by the mapping $(r,\theta) \in (0,+\infty) \times \s ^{n-1} \to r \theta$,
where $\nu$ is the one-dimensional probability measure on $\rr_+$ corresponding to the radial part, that is, with Lebesgue density proportional to $r^{n-1} e^{-V(r)}$.
As mentioned earlier, such a measure is invariant for the Sturm-Liouville dynamics
\begin{equation}
\label{eq:SL}
\LL _\nu f (r) := f''(r) - \left( V'(r) - \frac{n-1}{r} \right) f'(r) , \quad r >0,
\end{equation}
with additional Neumann's condition at the boundary $r=0$. Let $f \in \CC _0 ^\infty (\rr^n)$ and denote $g$ the bivariate function defined on $(0,+\infty) \times \s ^{n-1}$ by  $g(r,\theta) := f(r\theta)$. First using the Poincar\'e inequality \eqref{eq:poincare1d} for the radial operator $-\LL _\nu$, we have
\begin{eqnarray}
\label{eq:poincare_nu}
\nonumber \int _{\rr^n} f ^2 d\mu & = & \int _{\s ^{n-1}} \left( \int_0 ^{+\infty} g(r,\theta) ^2 \nu (dr) \right) \sigma_{n-1}(d\theta) \\
\nonumber & \leq & \int _{\s ^{n-1}} \left( \int_0 ^{+\infty} g(r,\theta) \nu (dr) \right) ^2 \sigma_{n-1}(d\theta) \\
& & + \frac{1}{\lambda_1 (-\LL_\nu)} \int _{\s ^{n-1}} \int_0 ^{+\infty} \left \vert \partial_r g(r,\theta) \right \vert ^2 \nu (dr) \sigma_{n-1}(d\theta) \\
\nonumber & = & \int _{\s ^{n-1}} \left( \int_0 ^{+\infty} g(r,\theta) \nu (dr) \right) ^2 \sigma_{n-1}(d\theta) \\
\nonumber & & + \frac{1}{\lambda_1 (-\LL_\nu)} \int _{\s ^{n-1}} \int_0 ^{+\infty} \left \vert \theta \cdot \nabla f (r\theta) \right \vert ^2 \nu (dr) \sigma_{n-1}(d\theta).
\end{eqnarray}
Now if we set
$$
h(\theta) := \int_0 ^{+\infty} g(r,\theta) \nu (dr) ,
$$
then the integral of $h$ under the measure $\sigma_{n-1}$ is nothing but the integral of $f$ under $\mu$. Moreover the length of the intrinsic gradient on the sphere applied to $h$ gives
\begin{eqnarray}
\label{eq:CS}
\nonumber \Vert \nabla _{\s ^{n-1}} h (\theta) \Vert ^2 & = & \left \Vert \int_0 ^{+\infty} r \Pi _{\theta ^\bot} (\nabla f) (r\theta) \nu (dr) \right \Vert ^2 \\
& \leq & \int_0 ^{+\infty} r^2 \nu (dr) \int_0 ^{+\infty} \left \Vert \Pi _{\theta ^\bot} (\nabla f) (r\theta) \right \Vert ^2 \nu (dr) \\
\nonumber & = & \int_{\rr ^n} \Vert x\Vert ^2 \mu (dx) \int_0 ^{+\infty} \left \Vert \Pi _{\theta ^\bot} (\nabla f) (r\theta) \right \Vert ^2 \nu (dr),
\end{eqnarray}
where to obtain \eqref{eq:CS} we used Cauchy-Schwarz inequality, the notation $\Pi _{\theta ^\bot}$ standing for the orthogonal projection on $\theta ^\bot$.
Since the Poincar\'e inequality \eqref{eq:poincare} is satisfied for the spherical Laplacian $- \Delta_{\s ^{n-1}}$, we have
\begin{eqnarray*}
\int _{\s ^{n-1}} h(\theta) ^2 \sigma_{n-1}(d\theta) & \leq & \left( \int _{\rr ^n } f d\mu \right) ^2 + \frac{1}{\lambda_1 (-\Delta_{\s ^{n-1}})} \int _{\s ^{n-1}} \Vert \nabla _{\s ^{n-1}} h (\theta) \Vert ^2 \sigma_{n-1}(d\theta) \\
& \leq & \left( \int _{\rr ^n } f d\mu \right) ^2 + \frac{\int_{\rr ^n} \Vert x\Vert ^2 \mu (dx) }{\lambda_1 (-\Delta_{\s ^{n-1}})} \int _{\s ^{n-1}} \int_0 ^{+\infty} \left \Vert \Pi _{\theta ^\bot} (\nabla f) (r\theta) \right \Vert ^2 \nu (dr) \sigma_{n-1}(d\theta) .
\end{eqnarray*}
Together with the previous estimate we obtain
\begin{eqnarray*}
\Var _\mu (f) & \leq & \frac{1}{\lambda_1 (-\LL_\nu)} \int _{\s ^{n-1}} \int_0 ^{+\infty} \left \vert \theta \cdot \nabla f (r\theta) \right \vert ^2 \nu (dr) \sigma_{n-1}(d\theta) \\
& & + \frac{\int_{\rr ^n} \Vert x\Vert ^2 \mu (dx) }{\lambda_1 (-\Delta_{\s ^{n-1}})} \int _{\s ^{n-1}} \int_0 ^{+\infty} \left \Vert \Pi _{\theta ^\bot} (\nabla f) (r\theta) \right \Vert ^2 \nu (dr) \sigma_{n-1}(d\theta) \\
& \leq & \max \left \{ \frac{1}{\lambda_1 (-\LL_\nu)} , \frac{\int_{\rr ^n} \Vert x\Vert ^2 \mu (dx) }{\lambda_1 (-\Delta_{\s ^{n-1}})} \right \} \int _{\s ^{n-1}} \int_0 ^{+\infty} \left \Vert \nabla f (r\theta) \right \Vert ^2 \nu (dr) \sigma_{n-1}(d\theta) ,
\end{eqnarray*}
since by the Pythagorean theorem,
\[
\Vert \nabla f \Vert ^2 = \left \vert \theta \cdot \nabla f \right \vert ^2 + \left \Vert \Pi _{\theta ^\bot} (\nabla f) \right \Vert ^2 .
\]
Therefore we get the inequality
\[
\lambda_1 (-\LL_\mu) \geq \min \left \{ \lambda_1 (-\LL_\nu) , \frac{\lambda_1 (-\Delta_{\s ^{n-1}})}{\int_{\rr ^n} \Vert x\Vert ^2 \mu (dx) }\right \} ,
\]
and finally using that $\lambda_1 (-\Delta_{\s ^{n-1}}) = n-1$ entails the desired lower bound. The proof is achieved.
\nprf

\subsection{Proof of Theorem \ref{theo:main}}
Once we have in mind the two previous lemmas, we can turn to the (brief) proof of our main result, Theorem~\ref{theo:main}.
\bprf[Proof of Theorem \ref{theo:main}]
As mentioned previously, the conclusion is a straightforward consequence of Lemmas \ref{lemme:bj} and \ref{lemme:gen-spherically} since we have
\[
\lambda_1 (-\LL_\nu) \geq \frac{n-1}{\int_{0}^{+\infty} r^2 \nu(dr)} = \frac{n-1}{\int_{\rr ^n} \Vert x\Vert ^2 \mu (dx) }.
\]
\nprf

\section{Examples}
\label{sect:ex}
In order to illustrate our Theorem \ref{theo:main}, let us revisit some classical models for which we are able to estimate the associated spectral gap.

\subsection{Uniform measure on the Euclidean unit ball}
First we consider the basic example of the Euclidean unit ball $B_n$ in $\rr^n$ ($n \geq 2$). In this case we have $V \equiv 0$ and $\mu$ is the normalized Lebesgue measure on $B_n$, i.e., its Lebesgue density on $B_n$ is given by $1/\omega_n$, with $\omega _n$ the volume of $B_n$,
$$
\omega _n = \frac{\pi ^{n/2}}{\Gamma \left( \frac{n}{2} +1 \right)},
$$
and where $\Gamma$ is the standard Gamma function
$$
\Gamma (z) := \int_0 ^{+\infty} x^{z-1} e^{-x} dx, \quad z>0.
$$
In terms of dynamics, the corresponding Markov process is the standard Brownian motion in $B_n$ reflected at the boundary (encoding Neumann's boundary conditions for the Euclidean Laplacian $\Delta$).
Then we obtain the following corollary of Theorem~\ref{theo:main}.
\begin{cor}[Euclidean unit ball]
Let $\lambda_1 (-\Delta)$ be the Neumann spectral gap associated to the Laplacian $\Delta$ and the normalized Lebesgue measure on the Euclidean unit ball $B_n$ in $\rr^n$ ($n\geq 2$). Then we have
$$
\frac{(n-1)(n+2)}{n} \leq \lambda_1 (-\Delta) \leq n+2.
$$
In particular $\lambda_1 (-\Delta) \approx n$ as the dimension $n$ is large.
\end{cor}

\subsection{Exponential power distributions}
Let us turn to another example containing the last one as a limiting case. We consider the potential $V(r) = r^\alpha /\alpha$ on $\rr_+$ which is convex provided $\alpha \geq 1$, so that the associated measure of interest $\mu$ has Lebesgue density on $\rr^n$ ($n\geq 2$) proportional to $e^{- \Vert x\Vert ^\alpha /\alpha}$. Although the terminology is not entirely stabilized in the literature, such probability measures are sometimes called exponential power distributions. The underlying dynamics is then the following:
\begin{equation}
\label{eq:exp_power}
\LL _\mu f (x) = \Delta f (x) - \Vert x \Vert ^{\alpha -2} x \cdot \nabla f (x) , \quad x\in \rr ^n .
\end{equation}
Among all these models, some particular cases reveal to be very interesting: \vspace{0.1cm}

$\circ$ When $\alpha = 2$, the associated Markov process is the Ornstein-Uhlenbeck process and $\mu$ is the centered Gaussian probability measure on $\rr^n$ with the identity matrix $I$ as covariance matrix. \vspace{0.1cm}

$\circ$ When $\alpha = 1$, the measure $\mu$ is an exponential-type distribution. In other words it is comparable to the product measure for which all the coordinates are i.i.d. and exponentially distributed with common parameter 1, since the $\ell ^1$ and $\ell ^2$ norms are equivalent in $\rr^n$. \vspace{0.1cm}

$\circ$ When $\alpha$ tends to infinity, we recover the uniform measure on the Euclidean unit ball $B_n$ studied previously. \vspace{0.1cm}

We still use the abuse of notation $V(x) = V(\Vert x\Vert )$. Then the Hessian matrix of $V$ at point $x \in \rr ^n $ is
$$
\mbox{Hess } V (x) = \Vert x\Vert ^{\alpha -2} I + (\alpha -2) \Vert x\Vert ^{\alpha -4} x x^T ,
$$
where $x^T$ stands for the transpose of the column vector $x$. The eigenvalues are $(\alpha -1) \Vert x\Vert ^{\alpha -2}$ with $x$ as corresponding eigenvector, and also $\Vert x\Vert ^{\alpha -2}$ associated to the eigenvectors orthogonal to $x$. Therefore we observe that Bakry-\'Emery criterion, which is reduced in our setting to the existence of some positive number $K$, independent from $x$, such that
$$
\mbox{Hess } V \geq K I,
$$
where the inequality has to be understood in the sense of positive matrices, cannot be applied except in the Gaussian case $\alpha = 2$ (for which it is sharp since we have the well-known result $\lambda_1 (-\LL_\mu) = 1$). Indeed, if $\alpha \in [1,2)$ then the smallest eigenvalue is $(\alpha -1) \Vert x\Vert ^{\alpha -2}$ which tends to 0 at infinity, whereas in the case $\alpha >2$, the smallest eigenvalue is now $\Vert x\Vert ^{\alpha -2}$ which vanishes at the origin. \vspace{0.1cm}

Coming back to the consequences of Theorem \ref{theo:main}, we obtain for this model the following result.
\begin{cor}[Exponential power distributions]
Let $\lambda_1 (-\LL_\mu)$ be the spectral gap of the operator $-\LL_\mu$ given in \eqref{eq:exp_power}. Then we have in any dimension $n \geq 2$,
\beq
\label{power-law1}
\frac{(n-1) \Gamma \left(\frac{n}{\alpha}\right)}{\alpha ^{2/\alpha} \Gamma \left( \frac{n+2}{\alpha}\right)} \leq \lambda_1 (-\LL_\mu) \leq \frac{n \Gamma \left(\frac{n}{\alpha}\right)}{\alpha ^{2/\alpha} \Gamma \left(\frac{n+2}{\alpha}\right)}.
\nneq
More explicitly, the two following estimates hold:
\beq
\label{power-law2}
\frac{n-1}{n+1} \times n^{1-2/\alpha} \leq \lambda_1(-\LL_{\mu}) \leq \frac{n+2}{n} \times n^{1-2/\alpha} .
\nneq
In particular we have $\lambda_1(-\LL_{\mu}) \approx n^{1-2/\alpha}$ as the dimension $n$ is large.
\end{cor}

\bprf
To obtain the inequalities \eqref{power-law1}, we use Theorem \ref{theo:main} and compute both upper and lower bounds. To establish the two last estimates, our idea is to invoke the log-convexity of the Gamma function. More precisely, if $a,b$ are two parameters such that $a>0$ and $b\in [0,1]$, then the the log-convexity of $\Gamma$ entails the two inequalities:
$$
\Gamma (a+b) \leq \Gamma (a) ^{1-b} \, \Gamma (a+1) ^b = a^b \, \Gamma (a) ,
$$
and
$$
\Gamma (a+1) \leq \Gamma (a+b) ^{b} \, \Gamma (a+b+1) ^{1-b} = (a+b)^{1-b} \, \Gamma (a+b) ,
$$
which rewrite in a condensed form as
\begin{equation}
\label{eq:gamma1}
1 \leq \frac{\Gamma (a) \, a^b}{\Gamma (a+b)} \leq \left( \frac{a+b}{a} \right) ^{1-b} .
\end{equation}
Moreover, if $b\in [1,2]$, then
$$
\frac{\Gamma (a) \, a^b}{\Gamma (a+b)} = \frac{a}{a+b-1} \times \frac{\Gamma (a) \, a^{b-1}}{\Gamma (a+b-1)} ,
$$
and since $b-1 \in [0,1]$, we can use the inequalities \eqref{eq:gamma1} to get
\begin{equation}
\label{eq:gamma2}
\frac{a}{a+b-1} \leq \frac{\Gamma (a) \, a^b}{\Gamma (a+b)} \leq \left( \frac{a+b-1}{a} \right) ^{2-b} .
\end{equation}
Now we use the previous estimates with $a = n/ \alpha >0$ and $b = 2/ \alpha \in (0,2]$. Then two cases occur: \vspace{0.1cm}

$\circ$ on the one hand, if $\alpha \geq 2$, that is $b \in (0,1]$, then the inequalities \eqref{eq:gamma1} entail the estimates
$$
1 \leq \frac{\Gamma \left( \frac{n}{\alpha} \right) \left( \frac{n}{\alpha}\right) ^{2/\alpha}}{\Gamma \left( \frac{n}{\alpha} + \frac{2}{\alpha}\right)} \leq \left( \frac{n+2}{n}\right) ^{1- 2/\alpha} \leq \frac{n+2}{n} .
$$

$\circ$ on the other hand, if $\alpha \in [1,2]$, i.e., $b\in [1,2]$, then we obtain from \eqref{eq:gamma2} the inequalities

$$
\frac{n}{n+2 - \alpha} \leq \frac{\Gamma \left( \frac{n}{\alpha} \right) \left( \frac{n}{\alpha}\right) ^{2/\alpha}}{\Gamma \left( \frac{n}{\alpha} + \frac{2}{\alpha}\right)} \leq \left( \frac{n+2 - \alpha }{n}\right) ^{2- 2/\alpha} .
$$
Finally, since in both cases the lower and upper bounds are respectively greater than or equal to $n/(n+1)$ and less than or equal to $(n+2)/n$, the desired estimates \eqref{power-law2} hold.
\nprf
As expected, the spectral gap is dimension-free only in the Gaussian case, i.e., $\alpha = 2$, reflecting that $\mu$ is a product measure.  When $\alpha = 1$ and the measure $\mu$ is of exponential-type, we obtain from \eqref{power-law1} the nice estimates
$$
\frac{n-1}{n(n+1)} \leq \lambda_1 (-\LL_\mu) \leq \frac{1}{n+1} .
$$

\subsection{Spectral gap of the one-dimensional radial part}
Let us briefly discuss the spectral gap $\lambda_1 (-\LL_\nu)$ of the associated one-dimensional radial operator $-\LL _\nu$. Recall that this generator is defined in \eqref{eq:SL} and for exponential power distributions, we have
$$
\LL _\nu f (r) := f''(r) - \left( r^{\alpha-1} - \frac{n-1}{r} \right) f'(r) , \quad r >0,
$$
with additional Neumann's condition at the boundary $r=0$. In contrast to the Cauchy-like case which will be analyzed in the next section, the exact value of the spectral gap does not seem to be accessible. However we are able to give its behaviour in terms of the underlying dimension $n$. Some simple computations in the inequality \eqref{eq:bj_radial} together with the choice of the function $f(r ) = r^\alpha$ in the variational inequality of the spectral gap show that $\lambda_1 (-\LL_\nu)$ is of order $n^{1-2/\alpha}$ for large $n$, as in the case of the full dynamics. For the exponential distribution, that is $\alpha =1$, such a result might be recovered easily. Indeed observing that the measure $\nu$, which is the Gamma distribution with $n$ degrees of freedom, is the $n$-fold convolution product of the standard exponential distribution, we obtain immediately that $\lambda_1 (-\LL_\nu)$ is at least $1/n$ times the spectral gap for the exponential measure, which is $1/4$. On the other hand, a curious phenomenon appears in the limiting case $\alpha \to +\infty$. In this situation, the measure $\nu$ converges weakly to the probability measure with density $n r^{n-1}$ on $(0,1)$ and the (Neumann) generator of interest is
$$
\LL _\nu f (r) := f''(r) + \frac{n-1}{r} \, f'(r) , \quad r \in (0,1).
$$
Then the inequality \eqref{eq:bj_radial} gives a lower bound of order $n$ on the spectral gap. However, taking the function $f' (r) = r^{(n-1)/2}$ in the inequality \eqref{e:chen}, one obtains
$$
V_f (r) = - \frac{(\LL _\nu f)'(r) }{f'(r) } = \frac{(n-1)^2}{4r^2} + \frac{n-1}{2r^2} \geq \frac{n^2 - 1}{4}, \quad r\in(0,1) ,
$$
and together with the choice of the identity function in the variational formula of the spectral gap show that the order of magnitude is $n^2$, which is different from the one of the full dynamics.

\section{Beyond the log-concave case}
\label{sect:beyond} The approach emphasized in the proof of Theorem~\ref{theo:main} allows us to obtain interesting results beyond the case of log-concave measures,
at the price of potentially modifying the energy term in the Poincar\'e inequality, leading henceforth to the so-called weighted Poincar\'e inequalities.
Let us observe how this phenomenon appears in the general setting before turning to two interesting examples: the first one is related to heavy-tailed measures, the generalized Cauchy distribution, whereas the second one revisits the Gaussian case.

\subsection{General setting}
As in the previous sections, the main protagonist of this part is still the probability measure $\mu$ on $\rr ^n$ ($n \geq 2$) with Lebesgue density proportional to $e^{-V (\Vert x\Vert)}$, but the notable difference relies in the smooth potential $V : \rr _+ \to \rr$ which is no longer convex.
If we consider as before the generator acting on smooth enough functions $f$,
$$
\LL _\mu f( x) := \Delta f (x) - \frac{V' (\Vert x\Vert )}{\Vert x\Vert } \, x  \cdot  \nabla f( x) , \quad x\in \rr^n ,
$$
then it may happen that $\lambda_1 (- \LL_\mu) = 0$, i.e., there is no spectral gap. In other words, there is no constant $\lambda >0$ such that the Poincar\'e inequality \eqref{eq:poincare} is satisfied. Among the potential reasons for which such a functional inequality can be false, a classical one is that the Euclidean norm might not be exponentially integrable with respect to $\mu$. However the story becomes different if the energy term is modified in a convenient way, that is to say, if we introduce a multiplicative (unbounded) weight $\sigma^2$ in the integral of the right-hand-side of \eqref{eq:poincare} to reach an inequality of the type
\begin{equation}
\label{eq:weigh_poinc}
\lambda \Var _\mu (h) \leq \int_{\rr^n} \sigma^2 \Vert \nabla h \Vert ^2 d\mu .
\end{equation}
Now, from our Markovian point of view, the natural question is the following: can we express the optimal constant in the latter weighted Poincar\'e inequality as the spectral gap of some alternative Markovian dynamics ? Actually the answer to this question is affirmative as we can see from now on. We introduce a new Markovian operator $\LL _\mu ^{\sigma}$ acting on functions $h\in \CC_0^\infty(\R^n)$ as follows,
\begin{equation}
\label{e:Lmusigma}
\LL ^\sigma _\mu h( x) = \sigma^2 (x) \Delta h (x) + \left( \nabla (\sigma^2) (x) - \sigma ^2(x) \frac{V' (\Vert x\Vert )}{\Vert x\Vert } x \right)  \cdot \nabla h( x) , \quad x\in \rr^n .
\end{equation}
Then the operator $\LL ^\sigma _\mu$ is symmetric and non-positive on $\CC_0^\infty(\R^n)$ with respect to the same measure probability $\mu$. Moreover, the main point resides in the fact that the weighted Poincar\'e inequality (\ref{eq:weigh_poinc}) rewrites by integration by parts as the following inequality: for smooth enough functions $h$,
$$
\lambda \Var _\mu (h) \leq \int _{\rr^n} h (-\LL ^\sigma _\mu h) d\mu ,
$$
henceforth leading to the study of the spectral gap $\lambda_1 (- \LL _\mu ^\sigma)$. \vspace{0.1cm}

Once we have realized that weighted Poincar\'e inequalities can be reduced to classical Poincar\'e inequalities for an alternative dynamics, the routine is more or less the same as in the previous section. First, let us introduce the following set of hypothesis on the weight function $\sigma$ that we will refer to as assumption $\AA$ in the sequel: \vspace{0.1cm}

$\circ$ Smoothness: the function $\sigma$ belongs to $\CC ^\infty (\R ^n)$. \vspace{0.1cm}

$\circ$ Ellipticity: it means that $\sigma$ is non-degenerate, i.e. $\sigma (x) > 0$ for every $x\in \rr ^n$. \vspace{0.1cm}

$\circ$ Completeness: the metric space $(\rr ^n, d_\sigma )$ is complete, where $d_\sigma $ is the distance
$$
d_\sigma (x,y) := \sup \, \left \{ \vert f(x) - f(y) \vert : f\in \CC ^\infty (\rr ^n) , \, \sigma \Vert \nabla f\Vert \leq 1 \right\} , \quad x,y\in \rr ^n .
$$
Then the operator $(\LL _\mu ^\sigma , \CC _0 ^\infty (\rr ^n))$ is essentially self-adjoint in $L^2(\mu)$. In this context, the Dirichlet form reads as
\[
 \EE_\mu^\sigma(f_1,f_2):= \int_{\rr ^n} \sigma ^2 \nabla f_1  \cdot  \nabla f_2 d\mu ,
\]
with domain $\DD (\EE _\mu ^\sigma)$ given by the weighted Sobolev space
\[
W_\sigma ^{1,2}(\mu) = \left\{ f\in L^2(\mu) : \sigma \Vert \nabla f \Vert \in  L^2(\mu) \right\} .
\]
As before, the self-adjointness allows us to check the Poincar\'e inequality \eqref{eq:weigh_poinc} only for functions $h\in \CC_0^\infty(\rr ^n)$. \vspace{0.1cm}

The second main object we have to deal with is the associated one-dimensional radial operator. To do so, we need a last assumption on the weight $\sigma$ which is the following: the function $\sigma$ is radial, i.e., it depends only on the length of the space variable. Recall that for radial functions, we use the same notation for the functions defined on $\rr^n$ and $(0,+\infty)$. Then the associated (Neumann) generator $\LL ^\sigma _\nu$ acts on smooth real-valued functions $g$ on $(0,+\infty)$ as
\begin{equation}
\label{e:Lnusigma}
\LL ^\sigma _\nu g (r) := \sigma ^2 (r) g''(r) + b(r)  g'( r) , \quad r >0,
\end{equation}
the drift $b$ being given by
$$
b(r) := (\sigma^2)' (r) - \sigma ^2 (r) \left( V'(r) - \frac{n-1}{r} \right) , \quad r >0.
$$

In contrast to the log-concave case studied before, there is no universal analogue of the one-dimensional estimate \eqref{eq:bj_radial}, so that we are not able to state an elegant result somewhat similar to Theorem~\ref{theo:main}. However in the spirit of Lemma~\ref{lemme:gen-spherically}, the spectral comparison between the full dynamics and its one-dimensional radial part is still possible and this is the content of the next result.
\bthm
\label{theo:gen-spherically-sigma}
Let $\mu$ be a spherically symmetric probability measure on $\R^n$ ($n \geq 2$) with Lebesgue density proportional to $e^{-V(\Vert x\Vert )}$ and let $\sigma$ be some radial function on $\rr ^n$ satisfying Assumption $\AA$. Let $\nu$ be the probability measure on $(0,+\infty)$ whose Lebesgue density is proportional to $r^{n-1} e^{- V(r)}$. Then provided the various integrals below make sense, the spectral gaps $\lambda_1 (-\LL_\mu)$ and $\lambda_1 (-\LL_\nu)$ of the operators $-\LL_\mu ^\sigma$ and $-\LL_\nu ^\sigma$ defined in \eqref{e:Lmusigma} and \eqref{e:Lnusigma} respectively, satisfy
\begin{equation}\label{e:min-maj-lambda_1-sigma}
\min \left \{ \lambda_1 (-\LL_\nu^\sigma) , \frac{n-1}{\int_{\rr ^n} \frac{\Vert x\Vert ^2 }{\sigma ^2(x) } \mu (dx)  }\right \}
\leq \lambda_1 (-\LL_\mu^\sigma)
\leq \min \left \{ \lambda_1 (-\LL_\nu^\sigma) , \frac{n \int_{\rr ^n} \sigma^2(x)  \mu (dx)} {\int_{\rr ^n} \Vert x\Vert ^2  \mu (dx)} \right \} .
\end{equation}
\nthm
\bprf
The proof is similar to the one of Lemma~\ref{lemme:gen-spherically}, with the following slight changes. First, the upper bound is obtained by choosing on the one hand radial functions $f$ in the variational characterization of the spectral gap $\lambda_1 (-\LL_\mu^\sigma)$ and on the other hand the linear function $h(x) = \sum_{k=1} ^n x_k$ in the weighted Poincar\'e inequality \eqref{eq:weigh_poinc}. \vspace{0.1cm}

For the lower bound, we follow the same computations as in the proof of Lemma~\ref{lemme:gen-spherically}, except in the Cauchy-Schwarz inequalities  \eqref{eq:poincare_nu}, \eqref{eq:CS} for which we exploit the presence of the weight $\sigma$, using $r = r /\sigma(r) \times \sigma(r)$. We arrive at the end at the estimate
\[
\lambda_1 (-\LL ^\sigma _\mu) \geq \min \left \{ \lambda_1 (-\LL ^\sigma _\nu) , \frac{\lambda_1 (-\Delta_{\sigma_{n-1}})}{\int_{\rr ^n} \frac{\Vert x\Vert ^2 }{\sigma ^2 (x) } \mu (dx) }\right \} ,
\]
which is the desired lower bound. The proof is now complete.
\nprf

Now, it is time to turn to the two examples we have in mind, which illustrate Theorem~\ref{theo:gen-spherically-sigma}: the heavy-tailed situation through generalized Cauchy distributions and then a new look at the Gaussian case.

\subsection{Generalized Cauchy distributions} First, let us focus our attention on the generalized Cauchy distribution. Let $\mu$ be the probability measure with Lebesgue density on $\rr^n$ ($n\geq 2$) proportional to $(1+\Vert x\Vert ^2) ^{-\beta}$, where $\beta >n/2$. In terms of the potential $V$, this density rewrites as $e^{-V(\Vert x\Vert )} /Z$ with
$$
V(r) := \beta \log (1+r^2), \quad r \geq 0,
$$
and where $Z$ stands for the normalization constant, that is
$$
Z := \frac{n \omega _n \, \Gamma (n/2) \, \Gamma(\beta - n/2)}{2 \Gamma(\beta)} .
$$
The measure $\mu$ is called the generalized Cauchy distribution of parameter $\beta$ on $\R^n$. \vspace{0.1cm}

Recently, Bobkov and Ledoux \cite{bob_ledoux2} used Brascamp-Lieb inequalities to establish, under the additional constraint $\beta \geq n$, the weighted Poincar\'e inequality \eqref{eq:weigh_poinc} with the radial weight
$$
\sigma^2 (x) := 1+ \Vert x\Vert ^2 ,
$$
and the constant $\lambda$ they obtained is
$$
\lambda = \frac{2(\beta -1)}{\left(\sqrt{1+\frac{2}{\beta -1}} + \sqrt{\frac{2}{\beta -1}} \right) ^2} .
$$
They also point out that a similar inequality might be established in the region $n/2 < \beta <n$, the constant $\lambda$ depending essentially on the dimension $n$. Since the weight $\sigma$ satisfies our Assumption $\AA$ and the generator $\LL_\mu ^\sigma$ rewrites as
\begin{equation}
\label{eq:sigma_mu_cauchy}
\LL_\mu ^\sigma h( x) = \left( 1+ \Vert x\Vert ^2 \right) \Delta h (x) + 2(1-\beta) x  \cdot \nabla h (x) , \quad x \in \rr ^n ,
\end{equation}
then their result means, in our Markovian language, that we have the spectral gap inequality:
$$
\lambda_1 (-\LL ^\sigma _\mu) \geq \frac{2(\beta -1)}{\left(\sqrt{1+\frac{2}{\beta -1}} + \sqrt{\frac{2}{\beta -1}} \right) ^2}.
$$
With this choice of weight function $\sigma$, we obtain after some computations involving the beta functions,
\[
\int_{\rr ^n} \Vert x\Vert ^2 \mu (dx) = \int_0^{+\infty} r^2 \, \nu(dr) = \frac{n}{2\beta - 2 - n},
\]
when $\beta > n/2 +1$ and also,
\[
\int_{\rr ^n} \frac{\Vert x\Vert ^2 }{\sigma ^2(x) } \mu (dx) = \int_0^{+\infty} \frac{r^2}{1+r^2}   \, \nu(dr) = \frac{n}{2\beta} ,
\]
so that applying Theorem~\ref{theo:gen-spherically-sigma} for this model yields, at least in the region $\beta > n/2 +1$,
\begin{equation}
\label{eq:thm_cauchy}
\min \left \{ \lambda_1 (-\LL_\nu^\sigma) , \frac{2\beta (n-1)}{n} \right \}
\leq \lambda_1 (-\LL_\mu^\sigma) \leq \min \left \{ \lambda_1 (-\LL_\nu^\sigma) , 2(\beta -1) \right \} .
\end{equation}
Recall that the associated one-dimensional (Neumann) radial generator $\LL_\nu$ is given by
\begin{equation}
\label{eq:sigma_nu_cauchy}
\LL_\nu ^\sigma g( r) = \left( 1+ r^2 \right) g''(r) + \left( (n+1 - 2\beta)r + \frac{n-1}{r} \right) g'(r) , \quad r >0.
\end{equation}
Hence it remains to estimate the spectral gap $\lambda_1 (-\LL ^\sigma _\nu)$. Surprisingly, its exact value can be computed for all $\beta > n/2$. On the one hand, consider the function
$$
g(r) := r^2 - \frac{n}{2\beta - 2 -n} , \quad r> 0.
$$
Such a function is nothing but an increasing eigenfunction of the operator $-\LL ^\sigma _\nu$ when $\beta > n/2 + 2$ and thus the corresponding eigenvalue is the spectral gap, that is,
$$
\lambda_1 (-\LL ^\sigma _\nu) = 4 \left( \beta - \frac{n}{2} - 1 \right) .
$$
On the other hand, for $\beta$ in the region $n/2 < \beta \leq n/2+2$, we claim that
\beq
\label{keybeta}
\lambda_1(-\LL_{\nu}^{\sigma}) = \left( \beta-\frac{n}{2} \right) ^2.
\nneq
In particular the spectral gap is continuous at the critical value $\beta = n/2+2$. To get the upper bound in \eqref{keybeta}, we choose $f(r):=(1+r^2)^{\vep}$ with $0<2\vep<\beta-n/2$ so that by the variational formula of the spectral gap, we have
\begin{eqnarray*}
\lambda_1(-\LL_{\nu}^{\sigma})^{-1} & \ge & \frac{\Var_{\nu}(f)}{\EE_{\nu}^{\sigma}(f, f)} \\
& =& \frac{\beta-2\vep}{2n\vep^2}\bigg(1-\frac{\Gamma (\beta-\vep-\frac n2)^2 \, \Gamma(\beta) \, \Gamma(\beta-2\vep)}{\Gamma(\beta-\frac n2) \, \Gamma(\beta-\vep)^2 \, \Gamma(\beta-2\vep-\frac n2)}\bigg),
\end{eqnarray*}
and then we take the limit as $\vep\to \beta /2 - n/4 $. To get the lower bound in \eqref{keybeta}, we choose $f(r):=(1+r^2)^{\frac{\beta}2-\frac n4}$. Using the formula \eqref{e:chen} for the dynamics $\LL_\nu ^\sigma$, we get after careful calculations,
\begin{eqnarray*}
V_f(r) & = & \frac{(-\LL_{\nu}^{\sigma} f)'(r)}{f'(r)} \\
& = & \frac{\left(\beta -\frac{n}{2}\right)^2 r^2 + 2\beta + n - n\beta + \frac{n^2}{2}}{1+r^2} ,
\end{eqnarray*}
leading to the estimate
\begin{eqnarray*}
\lambda_1(-\LL_{\nu}^{\sigma}) & \geq & \inf_{r>0} V_f(r) \\
& = & \left(\beta -\frac{n}{2}\right)^2 .
\end{eqnarray*}
Therefore combining with \eqref{eq:thm_cauchy} entails the following result for the generalized Cauchy distribution.

\begin{cor}[Generalized Cauchy distribution, case $n\geq 3$]
\label{cor:cauchy3}
Let $\mu$ be the generalized Cauchy distribution of parameter $\beta$ on $\R^n$ ($n\geq 3$). Consider the radial weight $\sigma ^2 (x):= 1+ \Vert x\Vert ^2$ and let $\LL_\mu ^\sigma$ and $\LL_\nu ^\sigma$ be the operators defined by \eqref{eq:sigma_mu_cauchy} and \eqref{eq:sigma_nu_cauchy} respectively. Then the spectral gaps $\lambda_1 (-\LL_\mu ^\sigma)$ and $\lambda_1 (-\LL_\mu ^\sigma)$ satisfy: \vspace{0.1cm}

$\circ$ for $n/2< \beta \leq n/2+2 $, 
$$
\lambda_1 (-\LL ^\sigma _\mu)= \lambda_1 (-\LL ^\sigma _\nu)=\left(\beta-\frac{n}{2}\right)^2;
$$

$\circ$ for $2+n/2 < \beta \leq n(n+2)/(n+1)$,
$$
\lambda_1 (-\LL ^\sigma _\mu) = \lambda_1 (-\LL ^\sigma _\nu) = 4 \left(\beta - \frac{n}{2}-1 \right);
$$

$\circ$ for $n(n+2)/(n+1) < \beta \leq n+1$, 
$$
\frac{2 \beta (n-1)}{n} \leq \lambda_1 (-\LL ^\sigma _\mu) \leq   \lambda_1 (-\LL ^\sigma _\nu) = 4 \left(\beta - \frac{n}{2}-1 \right);
$$

$\circ$ for $\beta > n+1$, 
$$
\frac{2 \beta (n-1)}{n} \leq \lambda_1 (-\LL ^\sigma _\mu) \leq 2(\beta-1) <  \lambda_1 (-\LL ^\sigma _\nu) = 4 \left(\beta - \frac{n}{2}-1 \right).
$$
\end{cor}
Let us say some words about this result. \vspace{0.1cm}

$\circ$ We ignore if we still have
$$
\lambda_1 (-\LL ^\sigma _\mu) = 4 \left(\beta - \frac{n}{2}-1 \right) ,
$$
for $\beta$ in the region $(n(n+2)/(n+1), n+1]$. However when $\beta > n+1$, we know that it is no longer the case since the upper bound $2(\beta-1)$ is slightly better, meaning that the value $\beta = n+1$ is critical for the spectral comparison between the full and the one-dimensional radial dynamics. \vspace{0.1cm}

$\circ$ Our lower bound $2\beta (n-1)/n$ is not numerically comparable to that of Bobkov-Ledoux (it depends on the value of $\beta$ with respect to the dimension $n$), but both exhibit more or less the same behaviour as $\beta$ is large (they are proportional to $\beta$). \vspace{0.1cm}

$\circ$ The reader may wonder whether or not the weight function $\sigma$ is optimal. As a partial answer, let us investigate the family of weights $\sigma_a ^2 (r) := (1+ r^2) ^a$ for $a\in (0,1)$. Choosing as above the function $f(r):=(1+r^2)^{\vep}$ with $0<2\vep<\beta-n/2$, we get by the variational formula of the spectral gap,
\begin{eqnarray*}
\lambda_1(-\LL_{\nu}^{\sigma_a})^{-1} & \ge & \frac{\Var_{\nu}(f)}{\EE_{\nu}^{\sigma_a}(f, f)} \\
& =& \bigg(\frac{\Gamma(\beta-2\vep-\frac n2)}{\Gamma(\beta-2\vep)}-\frac{\Gamma (\beta-\vep-\frac n2)^2 \, \Gamma(\beta)}{\Gamma(\beta-\frac n2) \, \Gamma(\beta-\vep)^2}\bigg) \frac{\Gamma(\beta-2\vep+2-a)}{2n\vep^2\Gamma(\beta-2\vep-\frac n2+1-a)} ,
\end{eqnarray*}
and then passing through the limit as $\vep\to \beta /2 - n/4 $, the right-hand-side tends to $+\infty$ or, in other words, $\lambda_1(-\LL_{\nu}^{\sigma_a}) = 0$. Hence one deduces that there is no spectral gap for those family of weights when $a \in (0,1)$. \vspace{0.1cm}

$\circ$ In the multi-dimensional case, we observe a rather different phenomenon from the one-dimensional case. As we have already seen in the one-dimensional case, the spectral gap has an associated eigenfunction (when it exists) which is strictly monotone. Reciprocally, if there exists a strictly monotone eigenfunction, then the associated eigenvalue is nothing but the spectral gap. In the Cauchy setting, it is worth noticing that for $\beta$ in the region $( 2+n/2 ,  n(n+2)/(n+1)]$, the spectral gap $4(\beta - n/2 -1)$ is associated to the eigenfunction, $h(x) = \Vert x\Vert ^2 - n/(2\beta -2 - n)$, which does not have monotonic coordinates. Moreover, in the case $n/2+1 < \beta < n+1$, the eigenvalue $2(\beta-1)$ is strictly greater than the spectral gap whereas the underlying eigenfunction, the linear function $h(x) = \sum_{k=1} ^n x_k$, has monotonic coordinates. \vspace{0.1cm}

$\circ$ It is interesting to observe that there is no eigenfunction associated to the spectral gap of the operator $-\LL ^\sigma _\nu$ in the case $n/2 <\beta < n/2 + 2$. Indeed, using the formula \eqref{e:chen} for the function $f(r)=(1+r^2)^{\frac{\beta }2-\frac{n}{4}}$, one has
\[
V_f(r)=  \frac{(-\LL_{\nu}^{\sigma} f)'(r)}{f'(r)} > \left( \beta-\frac n2 \right) ^2.
\]
Then the Brascamp-Lieb inequality established in \cite{bj} gives for all $h\in  W_\sigma ^{1,2}(\nu)$,
\[
 \Var_\mu(h) \leq \int_0^{+\infty} \frac{(h' )^2}{V_f } \sigma^2 d\nu < \frac{1}{\left(\beta-\frac{n}{2}\right)^2}\int_0^{+\infty} (h' ) ^2 \sigma^2 d\nu .
\]
This strict inequality prevents the existence of an eigenfunction with eigenvalue $ \left(\beta- n/2 \right)^2$. \vspace{0.1cm}

When $n=2,$ the above argument for the estimates of $\lambda_1(-\LL_{\mu}^{\sigma})$ and $\lambda_1(-\LL_{\nu}^{\sigma})$ is still valid, except the comparison between $\lambda_1(-\LL_{\nu}^{\sigma})$ and $2\beta (n-1)/n$. We thus obtain the following result.
\begin{cor}[Generalized Cauchy distribution, case $n = 2$]  When $n=2$, with the same notation as in Corollary \ref{cor:cauchy3}, we have: \vspace{0.1cm}

$\circ$ for $1<\beta\le \frac{3+\sqrt{5}}{2},$  \; $\lambda_1(-\LL_{\mu}^{\sigma})=\lambda_1(-\LL_{\nu}^{\sigma})=(\beta-1)^2;$\\

$\circ$ for $\frac{3+\sqrt{5}}{2}<\beta\le 3,$ \; $\beta\le \lambda_1(-\LL_{\mu}^{\sigma}) \le \lambda_1(-\LL_{\nu}^{\sigma})= (\beta-1)^2;$\\

$\circ$ for $\beta>3,$  \; $\beta\le\lambda_1(-\LL_{\mu}^{\sigma})\le 2(\beta-1)< \lambda_1(-\LL_{\nu}^{\sigma}) =4(\beta-2).$
\end{cor}

\subsection{The Gaussian case revisited}
As announced, we turn to our second example of interest, the case of weighted Poincar\'e inequalities for the Gaussian distribution. Let $\mu$ be the centered Gaussian probability measure on $\rr^n$ ($n\geq 2$) with the identity as covariance matrix. Then the question is the following: does the technology emphasized above allow us to obtain a weighted Poincar\'e inequality of the type \eqref{eq:weigh_poinc} for the Gaussian measure ? According to the previous discussion, our objective is to estimate the spectral gap $\lambda_1 (- \LL_\mu ^\sigma)$ of the modified dynamics $\LL_\mu ^\sigma$ defined in \eqref{e:Lmusigma}, which rewrites in the Gaussian context as
\begin{equation}
\label{e:sigma_mu_gauss}
\LL ^\sigma _\mu h( x) = \sigma^2 (x) \Delta h (x) + \left( \nabla (\sigma^2) (x) - \sigma (x) ^2 x \right)  \cdot \nabla h( x) , \quad x\in \rr^n ,
\end{equation}
the potential $V$ being now $V(r) := r^2/2$. \vspace{0.1cm}

Let us investigate first the case of a similar weight $\sigma$ as for the generalized Cauchy distributions above, that is
$$
\sigma^2 (x):= 1+ \Vert x\Vert ^2 ,
$$
so that the desired weighted Poincar\'e inequality is a weak version of the classical Poincar\'e inequality (different from the so-called weak Poincar\'e inequality). For the upper bound, we choose as above the linear function $h(x):= \sum_{k=1} ^n x_k$ in the weighted Poincar\'e inequality (\ref{eq:weigh_poinc}) to obtain the simple estimate,
\begin{eqnarray*}
\lambda_1 (- \LL_\mu ^\sigma) & \leq & n \times \frac{\int_{\rr^n} \left( 1 + \Vert x\Vert ^2\right) \mu (dx)}{\int_{\rr^n} \Vert x\Vert ^2 \mu (dx)} \\
& = & n+1.
\end{eqnarray*}
For the lower bound, observe that the lower bound in \eqref{e:min-maj-lambda_1-sigma} is still available in the weighted Gaussian case,
so that once again we are led to estimate the spectral gap $\lambda_1 (-\LL ^\sigma _\nu) $ of the one-dimensional radial operator $-\LL ^\sigma _\nu$. Since it exact value is not accessible, in contrast to the Cauchy case above, we need to revisit Lemma~\ref{lemme:bj} by taking in consideration the weight $\sigma$. We are thus led to the following result, which can be stated in a more general situation than in the Gaussian case (in particular, we decide to keep the notation involving $\sigma$ since it will also be useful when dealing with another choice of weight $\sigma$). Recall that in the Gaussian situation, the (Neumann) radial operator $\LL ^\sigma _\nu$ is given by
$$
\LL ^\sigma _\nu g (r) := \sigma ^2(r) g''(r) + b(r)  g'(r), \quad r >0,
$$
where the drift $b$ is defined by
$$
b(r) := (\sigma^2)' (r) - \sigma (r) ^2 \left( r - \frac{n-1}{r} \right) .
$$
\blem[Bonnefont-Joulin's lemma revisited \cite{bj}]
\label{lemme:bj_weight}
Assume that the function $\VV_\nu ^\sigma$ defined on $(0,+\infty)$ by
$$
\VV_\nu ^\sigma := \frac{\LL _\nu ^\sigma \sigma }{\sigma } - b' ,
$$
is positive. Then the spectral gap $\lambda_1 (-\LL_\nu ^\sigma)$ satisfies the inequality
\[
\lambda_1 (-\LL_\nu ^\sigma) \geq \frac{1}{\int_0 ^{+\infty} \frac{1}{\VV_\nu ^\sigma} d\nu } .
\]
\nlem

Coming back to our Gaussian setting, we have with this choice of weight $\sigma ^2 (r) = 1+ r^2$,
\begin{eqnarray*}
\VV_\nu ^\sigma (r) & = & \frac{4r^4 + 6r^2 +1}{1+r^2} + \frac{n-1}{r^2} \\
& \geq & 4r^2 ,
\end{eqnarray*}
so that applying Lemma~\ref{lemme:bj_weight} entails the following inequality, available for any $n \geq 3$,
\begin{eqnarray*}
\lambda_1 (-\LL _\nu ^\sigma) & \geq & \frac{\int_0 ^{+\infty} r^{n-1} e^{-r^2 /2} dr }{\int_0 ^{+\infty} \frac{r^{n-1}}{4r^2} e^{-r^2 /2} dr} \\
& = &  4(n-2) .
\end{eqnarray*}
For $n=2$, we have
$$
\VV_\nu ^\sigma (r) \geq 4r^2 + \frac{1}{r^2} \geq 4 .
$$
Hence by the lower bound in \eqref{e:min-maj-lambda_1-sigma}, we obtain for any $n\geq 2$ the estimate
\begin{eqnarray*}
\lambda_1 (-\LL ^\sigma _\mu) & \geq & \min \left \{ \lambda_1 (-\LL _\nu ^\sigma) , (n-1) \times \frac{ \int_0^{+\infty} r^{n-1} e^{-r^2 /2}dr}{\int_0 ^{+\infty} \frac{r^2}{ 1+r^2} r^{n-1} e^{-r^2 /2}dr }\right \} \\
& \geq & n-1,
\end{eqnarray*}
since $r^2 /(1+r^2) \leq 1$. Let us summarize our result in the weighted Gaussian case.
\begin{cor}[Weighted Gaussian case 1]
Let $\mu$ be the standard Gaussian distribution on $\rr^n$ ($n \geq 2$) and consider the radial weight $\sigma ^2 (x) := 1 + \Vert x \Vert ^2$. Let $\lambda_1 (-\LL_\mu ^\sigma)$ be the spectral gap of the operator $-\LL_\mu ^\sigma$ defined in \eqref{e:sigma_mu_gauss}. Then we have
$$
n-1 \leq \lambda_1 (-\LL ^\sigma _\mu) \leq n+1.
$$
In particular $\lambda_1 (-\LL ^\sigma _\mu)  \approx n$ in large dimension.
\end{cor}

A more challenging problem is to obtain a weighted Poincar\'e inequality for the Gaussian measure that is a reinforcement of the classical Poincar\'e inequality. Motivated by an open question raised in \cite{bob_ledoux2}, our second choice of radial weight will be
$$
\sigma ^2 (x) := \frac{1}{1+\Vert x\Vert ^2} ,
$$
which satisfies also Assumption $\AA$. For bounding from above the spectral gap $\lambda_1 (- \LL_\mu ^\sigma)$, we choose as before the linear function $h(x) = \sum_{k=1} ^n x_k$ in the weighted Poincar\'e inequality (\ref{eq:weigh_poinc}) to get the following upper bound available for any $n \geq 3$,
\begin{eqnarray*}
\lambda_1 (- \LL_\mu ^\sigma) & \leq & n \times \frac{\int_{\rr^n} \left( 1 + \Vert x\Vert ^2\right) ^{-1} \mu (dx)}{\int_{\rr^n} \Vert x\Vert ^2 \mu (dx)} \\
& \leq & n \times \frac{\int_{\rr^n} \Vert x\Vert ^{-2} \mu (dx)}{\int_{\rr^n} \Vert x\Vert ^2 \mu (dx)} \\
& = & \frac{1}{n-2}.
\end{eqnarray*}
For $n=2$ we use the trivial inequality $\left( 1 + \Vert x\Vert ^2\right) ^{-1} \leq 1$ to get $\lambda_1 (- \LL_\mu ^\sigma) \leq 1$. \vspace{0.1cm}

For the lower bound in the general case $n\geq 2$, we apply Lemma~\ref{lemme:bj_weight} in our setting with the positive potential
$$
\VV_\nu ^\sigma (r) = \frac{(2n-3)r^4 + (3n-1)r^2 +n-1}{\left( 1+r^2 \right) ^3 r^2} , \quad r>0,
$$
and using the inequality
$$
(2n-3)r^4 + (3n-1)r^2 +n-1 \geq (n-1)(1+r^2)^2 , \quad r >0,
$$
we obtain by the lower bound in \eqref{e:min-maj-lambda_1-sigma} the estimate
\begin{eqnarray*}
\lambda_1 (-\LL ^\sigma _\mu) & \geq & \min \left \{ \frac{\int_0^{+\infty} r^{n-1} e^{-r^2 /2}dr}{\int_0 ^{+\infty} \frac{\left( 1+r^2 \right) ^3 r^{n+1}}{(2n-3)r^4 + (3n-1)r^2 +n-1} e^{-r^2 /2}dr} , (n-1)  \times \frac{\int_0^{+\infty} r^{n-1} e^{-r^2 /2}dr}{\int_0 ^{+\infty} ( 1+r^2 ) r^{n+1} e^{-r^2 /2}dr }\right \} \\
& = & (n-1) \times \frac{ \int_0^{+\infty} r^{n-1} e^{-r^2 /2}dr}{\int_0 ^{+\infty} ( 1+r^2 ) r^{n+1} e^{-r^2 /2}dr } \\
& = & \frac{n-1}{n(n+3)}.
\end{eqnarray*}
Finally we get the following estimates, exhibiting a different behaviour with respect to the dimension.
\begin{cor}[Weighted Gaussian case 2]
Let $\mu$ be the standard Gaussian distribution on $\rr^n$ ($n \geq 2$) and consider the radial weight $\sigma ^2 (x) := 1/(1 + \Vert x \Vert ^2)$. Let $\lambda_1 (-\LL_\mu ^\sigma)$ be the spectral gap of the operator $-\LL_\mu ^\sigma$ defined in \eqref{e:sigma_mu_gauss}. Then we have
$$
\frac{n-1}{n(n+3)} \leq \lambda_1 (-\LL ^\sigma _\mu) \leq \frac{1}{n-2} \wedge 1 ,
$$
so that $ \lambda_1 (-\LL ^\sigma _\mu) \approx 1/n$ for large $n$.
\end{cor}

\end{document}